\documentclass{article}
\usepackage{amssymb,amsmath,graphicx,ucs,hyperref}
\usepackage[utf8x]{inputenc}

\newtheorem{theorem}{Theorem}

\title{Compact packings of space\\with three sizes of spheres}
\author{Thomas Fernique}
\date{}

\begin{document}

\maketitle

\begin{abstract}
A sphere packing of the three-dimensional Euclidean space  is compact if it has only tetrahedral holes, {\em i.e.}, any local maximum of the distance to the spheres is at equal distance to exactly four spheres.
This papers describes all the compact packings with spheres of three different sizes.
They are close-compact packings of unit spheres with holes filled in four different ways by smaller spheres.
\end{abstract}


\section{Introduction}
\label{sec:intro}

A {\bf packing} is a set of interior disjoint closed balls of $\mathbb{R}^3$.
Balls are usually called spheres in this context, and we shall also call them {\em beads} in this paper.
The {\bf contact graph} is the graph whose vertices are the sphere centers and the edges connect two vertices whenever the two spheres they are center of are tangent.
Such a packing is said to be {\bf compact} if its contact graph is a homogenoeous simplicial complex of dimension $3$.
In other words, it can be seen as a {\bf tiling} by the tetrahedra whose edges connect the centers of four pairwise tangent spheres, that is, a covering by such tetrahedra, with the intersection of two tetrahedra being a face, an edge or a vertex.\\

There is no compact packing with unit spheres because regular tetrahedra do not tile the space.
In \cite{Fer19}, we proved that the only compact packings with two sizes of spheres are obtained by filling the octahedral holes of a close-packing of unit spheres by spheres of size $r=\sqrt{2}-1$.
We prove here that the situation is quite similar with three size of spheres: the only compact packings are indeed obtained by filling holes of a close-packing of unit spheres:

\begin{theorem}
\label{th:main}
There are four pairs $(r,s)$ which allow compact packings with spheres of size $s<r<1$.
These compact packings are obtained from a close-compact packing of unit spheres as follows:
\begin{enumerate}
\item
Fill all the octahedral holes by a sphere of size $r$.
Fill tetrahedral holes by a sphere of size $s$.
$$
r=\sqrt{2}-1\simeq 0.41
\quad\textrm{and}\quad
s=\sqrt{\tfrac{3}{2}}-1\simeq 0.22
$$
\item
Fill all the octahedral holes by a sphere of size $r$.
Fill tetrahedral holes by a tetrahedron of four spheres of size $s$.
$$
r=\sqrt{2}-1\simeq 0.41
\quad\textrm{and}\quad
s=3-2\sqrt{2}\simeq 0.17
$$
\item
Fill all the octahedral holes by a sphere of size $r$.
Each octahedral holes yields four tetrahedra 111r.
Fill such new tetrahedral holes by a sphere of size $s$.
$$
r=\sqrt{2}-1\simeq 0.41
\quad\textrm{and}\quad
s=\tfrac{1}{2}\sqrt{9+6\sqrt{2}}-\sqrt{2}-\tfrac{1}{2}\simeq 0.18
$$
\item
Fill all the octahedral holes by a tetrahedron of four spheres of size $r$.
This creates, besides this tetrahedron rrrr, four tetrahedra 111r and four octahedra with a face rrr and a face 111.
Fill all these new octahedral holes by a sphere of size $s$.
$$
r=4-\sqrt{14}\simeq 0.26
\quad\textrm{and}\quad
s=6-3\sqrt{2}-\sqrt{42-28\sqrt{2}}\simeq 0.21
$$
\end{enumerate}
\end{theorem}

In all these compact packings, the largest spheres form a close-packing and the smaller spheres fill holes.
This is a bit disappointing, since we could have expect new structures which cannot be that easily derived from the usual close-packing.
Does exist a set of spheres which allow a compact packing where it is not the case?
In particular, do four sizes of spheres suffice?

\section{Geometry of tiles}
\label{sec:geometry}

Compact packings can be seen as tilings by tetrahedra whose edges connect the centers of four pairwise tangent balls.
There are, up to isometry, finitely many different tetrahedra (actually $15$), depending on the radii of the involved balls.
Since we shall further rely on the dihedral and solid angles of these tetrahedra, let us here show how to use spherical trigonometry to compute them.\\

Let $A$, $B$, $C$ and $D$ be the centers of pairwise tangent balls of respective radii $a$, $b$, $c$ and $d$.
All what follows holds up to a permutation of the vertices.
Denote by
\begin{itemize}
\item $\widehat{abc}$ the angle in $B$ in the face $ABC$;
\item $\widehat{ab/cd}$ the dihedral angle between the faces intersecting on the edge $AB$;
\item $\widehat{a/bcd}$ the solid angle in $A$.
\end{itemize}
The (planar) law of cosines first yields:
$$
\cos\widehat{abc}=\frac{(b+a)^2+(b+c)^2-(a+c)^2}{2(b+a)(b+c)}.
$$
The spherical law of cosines then yields:
$$
\cos \widehat{ab/cd}=\frac{\cos\widehat{dac}-\cos\widehat{dab}\times\cos\widehat{bac}}{\sin\widehat{dab}\times\sin\widehat{bac}}.
$$
The Girard's theorem finally yields:
$$
\widehat{a/bcd}=\widehat{ab/cd}+\widehat{ac/bd}+\widehat{ad/bc}-\pi.
$$
All the dihedral and solid angles can thus be derived from the radii.\\

Since $\widehat{abc}=\widehat{cba}$, we can assume $a\geq c$.
This yields $6$ choices for $ac$ and $3$ for $b$, hence at most $18$ different angles, actually $16$ because
$$
\widehat{111}=\widehat{rrr}=\widehat{sss}=\tfrac{\pi}{3}.
$$
Since $\widehat{ab/cd}$ is unchanged by permuting $a$ and $b$ or $c$ and $d$, we can assume $a\geq b$ and $c\geq d$.
There is thus at most $6\times 6$ different dihedral angles, actually $34$ because
$$
\widehat{11/11}=\widehat{rr/rr}=\widehat{ss/ss}=\arccos\tfrac{1}{3}.
$$
Since $\widehat{a/bcd}$ is unchanged by a permutation of $\{b,c,d\}$, we can assume $b\geq c\geq d$.
This yields $10$ choices for $abc$, hence at most $30$ different solid angles, actually $28$ because
$$
\widehat{1/111}=\widehat{r/rrr}=\widehat{s/sss}=\arccos\tfrac{23}{27}.
$$

\section{Necklaces}
\label{sec:necklaces}

If $B$ and $H$ are two adjacent balls of a compact packing, consider the sequence $B_1,\ldots,B_k$ of balls which are each adjacent to both $B$ and $H$, as well as to both $B_i$ is adjacent to $B_{i+1}$.
We call the sequence of the $B_i$'s a {\bf necklace}, seing $B$ as the body, $H$ as the head and the $B_i$'s as the beads.
In other words, the $B_i$'s are centered on the vertices of the tetrahedra (other than the centers of $B$ and $H$) which share the edge which connects the centers of $B$ and $H$.\\

A necklace is coded by the word over $\{1,r,s\}$ whose $i$-th letter gives the radius of the $i$-th bead.
Among the possible codings, we usually choose the lexicographically minimal one.
We speak about a bh-necklace when $B$ and $H$ have radius $b$ and $h$ (or conversely).
There are thus $6$ different types of necklaces.\\

\noindent The bh-necklace $i_1\cdots i_k$ yields the equation:
$$
\widehat{bh/i_1i_2}+\widehat{bh/i_2i_3}+\ldots+\widehat{bh/i_ki_1}=2\pi.
$$
The planar and spherical laws of cosines then yield an equation in $r$ and $s$.
It is however generally not algebraic, and we prefer to associate with such a necklace the system of polynomial equations defined by the following steps:
\begin{enumerate}
\item take the cosines of both sides of the above equation;
\item fully expand the left-hand side yields a polynomial in cosines and sines of dihedral angles:
replace the cosine and sine of the $i$-th angle by auxiliary variables $X_i$ and $Y_i$ to get the first polynomial equation of the system;
\item add the equations $X_i^2+Y_i^2=1$ to the system;
\item the spherical law of cosines gives $Y_i$ as a rational fraction of cosines and sines of angles of tetrahedron faces: replace the latter by auxiliary variables $Z_j$ and $T_j$ and multiply both sides of the equation by the lowest common denominator to get a new polynomial equation added to the system;
\item add the equations $Z_j^2+T_j^2=1$ to the system;
\item the planar law of cosines gives $Z_j$ as a rational fraction in $r$ and $s$: multiply both sides of the equation by the lowest common denominator to get a polynomial equation and add it to the system.
\end{enumerate}
The advantage is that systems of polynomial equations can be exactly solved (at least in theory) by modern computer algebra systems when they are $0$-dimensional, {\em i.e.}, admit only finitely many solutions. 
Because of Step~$1$, though, any solution of this system yields values of $r$ and $s$ which are only solutions for some integer $n$ of the equation
$$
\widehat{bh/i_1i_2}+\widehat{bh/i_2i_3}+\ldots+\widehat{bh/i_ki_1}=2n\pi.
$$
To find the solution of the original equation, we shall then check that $n=1$.
This can be easily done from the exact values of $r$ and $s$ by using interval arithmetic with an accuracy only greater than $2\pi$.\\

Note that only the first equation of this system (Step~$2$) is specific to the necklace.
All the other ones are generic and can be precomputed.
Since there are at most $34$ dihedral angles and $16$ angles, this yields at most $100$ such equations.

\section{Shells}
\label{sec:shells}

Given a ball $C$ in a compact packing, the {\bf shell} is the set of balls of this packing which are in contact with $C$.
The center of $C$ is called the center of the shell.\\

A shell is encoded by a labelled spherical triangulation: each sphere around $C$ corresponds to a vertex labelled with the radius of this sphere, and edges connect the vertices corresponding to adjacent spheres.
We speak about a $c$-shell when $C$ has radius $c$.\\

\noindent A $c$-shell represented by a triangulation $T$ of size $k$ yields the equation
$$
\sum_{UVW\in T}\widehat{c/uvw}=4\pi,
$$
where $u$, $v$ and $w$ are the labels of the vertices $U$, $V$ and $W$.
This can be rewritten via Girard's theorem
$$
\sum_{uvw\in T}\widehat{cu/vw}+\widehat{cv/uw}+\widehat{cw/uv}=(4+k)\pi.
$$
The same procedure as for the necklaces transforms this equation in a system of polynomial equations whose solutions yield the solutions of the original equation.

\section{Large-separated packings}
\label{sec:large_separated}

A compact packing without contact between small and medium spheres is said to be {\em large-seperated}.
We here prove that there are only two set of sizes which allow such packings.
We rely on intermediary results proven in \cite{Fer19} for $r<1$:
\begin{enumerate}
\item there are only $10$ possible 1r-necklaces made of spheres of sizes $1$ and $r$;
\item only one of these $10$ 1r-necklace can coexist in the same packing with an rr-necklace, namely the 1r-necklace 11rr can coexist with the 1r-necklace 11rr for $r=r_{tt}:=3-2\sqrt{2}\simeq 0.17$.
\end{enumerate}

\paragraph{r-shells.}
Consider a r-shell in such a packing.
It contains only spheres of size $r$ and $1$.
It cannot contains only spheres of size $r$ because the solid angle of a regular tetrahedron does not divide $4\pi$.
It thus contains a large sphere, hence a 1r-necklace.
This 1r-necklace contains only spheres of the shell, thus of sizes $r$ and $1$.

If $r\neq r_{tt}$, then there is no rr-necklace, so that the shell actually contains only spheres of size $1$.
The polyhedron which connects the center of these spheres is a triangulated regular polyhedron, hence either a tetrahedron, an octahedron or an icosahedron.
The respective values of $r$ are
$$
r_t:=\sqrt{\tfrac{3}{2}}-1\simeq 0.22,
\qquad
r_o:=\sqrt{2}-1\simeq 0.41,
\qquad
r_i=\sqrt{\tfrac{1+\sqrt{5}}{2}+2}-1\simeq 0.90.
$$

If $r=r_{tt}$, the rr-necklace 11rr is allowed (on no other rr-necklace).
The shell contains at least one r-sphere, otherwise we are in the former case which appeared to be incompatible with $r=r_{tt}$.
Consider an r-sphere in the shell.
It is surrounded (in the shell) by the necklace 11rr.
Each sphere in this necklace is itself surrounded by 11rr.
This show that the shell actually contains only $6$ spheres: three r-spheres and three 1-spheres.
The three r-spheres form, with the center of the shell, a tetrahedron of r-spheres.
The shell of each of the four r-spheres of this tetrahedron contains only 1- and r-spheres: the same reasoning as above yields the same shell.
This tetrahedron of four r-spheres is thus surrounded by a "cage" of 1-spheres.
Since each of these four r-spheres is surrounded by four 1-spheres and each of these surrounding 1-spheres is adjacent to three of the four r-spheres, this "cage" is made of four 1-spheres.
In other words, the r-spheres form a tetrahedron filling the hole in a tetrahedron of 1-spheres.

\paragraph{s-shells.}
When we considered an r-shell, only 1- and r-spheres were involved.
The same reasoning thus applies for s-shells, just with $s$ instead of $r$ everywhere.

\paragraph{Packings.}
Any large-seperated compact packing can be seen as a tiling by two polyhedra among the tetrahedron, the octahedron and the icosahedron (two different ones since $s\neq r$), with 1-spheres centered on the vertices of these polyhedra and r- and s-spheres filling holes of these polyhedra.
Such tilings are called {\em uniform honeycomb} and have been classified in \cite{Gru94}.
The only possible case are {\em tetrahedral-octahedral honeycombs}, where vertices of tiles are the center of a close-packing of spheres.
This yields the two first cases of Theorem~\ref{th:main}.

\section{Medium-separated packings}
\label{sec:medium_separated}

A compact packing without contact between small and large spheres is said to be {\em medium-seperated}.
We here prove that there is no such packing.

\paragraph{s-shells.}
Up to a scale factor $\tfrac{s}{r}$, the situation for s-shells is the same as for r-shells in large-sperated packings.
We can thus have either s-shells with r-spheres centered on the vertices of a tetrahedron, or an octahedron or an icosahedron (depending whther $\tfrac{s}{r}$ is equal to $r_t$, $r_o$ and $r_i$), or s-spheres forming tetrahedra which fill the holes between tetrahedra of r-spheres.

\paragraph{1-shells.}
For a 1-shell, the situation is slightly different.
It contains only spheres of size $1$ and $r$ and cannot contain only large spheres.
There is thus a 1r-necklace with spheres of sizes $r$ and $1$.
We shall again rely on intermediary results proven in \cite{Fer19} for $r<1$:
\begin{enumerate}
\item Any of the $10$ possible 1r-necklaces containing only spheres of sizes $1$ and $r$ contains at least one sphere of size $1$;
\item only $r=r_o$ allows a 1r- and a 11-necklace both containing only spheres of sizes $1$ and $r$.
\end{enumerate}
Thus $r=r_o$ and the only possible 1-shells are those appearing in a compact packing by two sizes of spheres \cite{Fer19} (they both contain $12$ large spheres and $6$ medium ones, respectively centered on the vertices of either a cubocatahedron and an octahedron or a triangular orthobicupola and a triangular prism).

\paragraph{Packings.}
The 1-spheres form a close-compact packings whose octahedral holes are filled by r-spheres (same proof as in \cite{Fer19}).
The s-spheres must then fill holes between "cages" of r-spheres.
Unfortunately there is no such "cage" (each r-sphere is itself in a octahedral holes of 1-spheres).
There is thus no medium-separated packings with three sizes of discs.

\section{The remaining packings}
\label{sec:remaining}

\paragraph{1s- and rs-necklaces.}
The remaining packings are those with both 1s- and rs-necklaces.
The point is that we can bound the number of 1s- and rs-necklaces (as well as ss-necklaces) independently of the value of $s$ (whereas 11-, 1r- and rr-necklaces can contain arbitrarily many s-beads for $s$ small enough).
More precisely, 1s- or rs-necklace contains three and five beads, and an exhaustive search yields $70$ different codings.
This thus yields $70^2=4900$ pairs of 1s- and rs-necklaces to consider.

\paragraph{Neighboor necklaces.}
In a packing, two necklaces which share the same body and whose heads are adjacent are said to be {\em neighboor}.
Two such necklaces have two beads in common, say a- and b-beads, and if their heads are x- and y-beads, then their respective codings contain factors ayb and axb.
We shall use this constraint to reduce the number of pair we have to consider.
Consider indeed one of the above $4900$ pairs of 1s- and rs-necklace.
If the 1s-necklace contains an r-bead, then any packing which contains this 1s-necklace must also contain a neighboor rs-necklace (and symmetrically if the rs-necklace contains an 1-bead).
In such a case, we can thus assume that the pair of necklaces which characterize such packings are neighboor and check that the necklaces satisfy the above constraint.
If not, we just rule out the pair (because there is another pair which will describe the same packings).
This is easily automated and enventually reduce to $1103$ pairs of 1s- and rs-necklaces.

\paragraph{Computations.}
We then proceed similarly to \cite{FHS18}:
\begin{enumerate}
\item
We associate with each necklace a polynomial equation in r and s as in  (with dihedral angles instead of angles and the spherical cosine law playing the role of the cosine law).
The mean degree over the $70$ possible 1s- or rs-necklaces is around $16$, with maximum $72$ for the necklace 1rrss (both in the 1s- and rs-case, though polynomials are different).
\item
We use the hidden variable method, {\em i.e.}, we compute the two resultants in r and s of the two polynomials associated with a pair of 1s- and rs-necklaces, compute exactly the roots, filter the roots $0<s<r<1$, then filter by checking with interval arithmetic on dihedral angles that they are compatible with the 1s- and rs-necklaces the polynomials have been associated with.
\item
If an s-bead appears in the 1s- or rs-necklace, then we search for a compatible ss-necklace (with interval arithmetic on dihedral angles).
If an r-bead appears in the 1s-necklace or a 1-bead appears in the rs-necklace, then we search for a compatible 1r-necklace.
If a 1-bead appears in the 1s-necklace, then we search for a compatible 11-necklace.
If an r-bead appears in the rs-necklace, then we search for a compatible rr-necklace.
\end{enumerate}
This filtering suffices to rule out all the pairs but $15$.
One is the pair 11ss/rrrr.
The $14$ other ones yield a polynomial system of dimension one.

\paragraph{The pair 11ss/rrrr.}
The associated equations yield the values
$$
r=\sqrt{2}-1\simeq 0.41
\qquad
s=3-2\sqrt{2}\simeq 0.17.
$$
These are the values $r_o$ and $r_{tt}$ which already appeared in Section~\ref{sec:large_separated}.
We shall show that we get the same packings, that is, the possible packings are actually large-separated.
The values of $r$ and $s$ allows to check by computation that 11ss (resp. rrrr) is the unique 1s-necklace (resp. rs-necklace).
They also allows to compute solid angles, hence to search for the number of triangles of each type which could form an s-shell.
The computation yields two cases:
\begin{enumerate}
\item $8$ triangles rrr;
\item $1$ triangle 111, $3$ triangles 11s, $3$ triangles 1ss and $1$ triangle sss.
\end{enumerate}
In the first case, the $s$-shell is an octahedron of r-beads.
In the second case, the s-shell is an octahedra with a face 111 and a face sss.
In this latter case, an s-bead of the sss face cannot have an octahedral shell of r-beads (because it is adjacent to an s-bead), so it has the same shell as the first s-bead.
This enforces the four s-beads to form a tetrahedra around which 1-beads form a larger tetrahedra.
Now, assume there is a compact packing with such beads (where each type of bead does appear).
Remove the tetrahedra of s-beads which fill tetrahedra of 1-beads (if any): the packing is still compact.
Now, the 1-shell contain only 1- and r-beads (because the shell of any remaining s-bead is an octahedron of r-beads, so no s-bead can be in contact with an 1-bead).
A computation shows that there are two 11-necklaces without s-bead, namely 111r1r and 11r11r, and only one 1r-necklace without s-bead, namely 1111.
These necklaces are exactly those we got in \cite{Fer19} while searching for compact packings with spheres of size $r=\sqrt{2}-1$ and $1$.
We there proved that the possible 1-shells are those appearing in a close-packing of 1-beads with each octahedral hole filled by an r-bead.
We also proved that this suffices to ensure that the whole packing is indeed a close-packing of 1-beads with each octahedral hole filled by an r-bead, without relying on r-shells (thus, even if r-shells which contain s-beads are possible, they cannot appear in a compact packing which contains a 1-bead).    
Now, if we put back in the tetrahedral holes of 1-beads the tetrahedra of s-beads that we removed, we fall back into the second case of Theorem~\ref{th:main}.

\paragraph{Dimension one: when there must be an ss-necklace.}
Out of the $14$ cases which yield a dimension $1$ system, four of them have an s-bead in the 1s-necklace or in the rs-necklace:
$$
rrs/1rs,
\qquad
1rs/11s,
\qquad
rss/1ss,
\qquad
rsrs/1s1s.
$$
Sizes of beads thus must allow an ss-necklace, which yields a third equation.
In each case, we consider the systems formed by adding the equation assoiated with one of the $70$ possible necklaces.
We check (as previously with resultant and interval airthmetic filtering) that none of these systems admits a solution, except one system in each of these four cases, which still has dimension one.
Namely, these are the four systems associated with the following 1s/rs/ss-necklaces:
$$
rrs/1rs/1rr,
\qquad
1rs/11s/11r,
\qquad
rss/1ss/1rs,
\quad
rsrs/1s1s/1r1r.
$$
Each of these four cases determine a unique s-shell.
These s-shells are, respectively the tetrahedra 1rrs, 11rs, 1rss and the octahedron with two beads of each type and no two adjacent identical beads.
The two first tetrahedra are impossible because the convex hull of the centers of the beads of the s-shell does not contains the center of the central s-bead (the centers of the central s-bead and the s-bead in the s-shell are on both side of the plane containing the centers of the 1- and r-beads), whereas this should but the case if the contact graph is a tiling by tetrahedra.
For the third tetrahedron, there must be three mutually adjacent s-beads and two mutually adjacent 1- and r-beads on both side of this triangle of s-beads: this is impossible because the 1- and r-beads are larger than the s-beads (thus cannot be adjacent in the hole between the three s-beads).
For the octahedron, its symmetry ensures that the three s-beads (the central one and the two in the shell itself) are aligned, say vertically.
The same argument then holds for the top-most s-bead: it is aligned with the underneath s-bead (the center of the former s-shell) and with a higher s-bead.
This latter s-bead is still adjacent to the two 1-beads and two r-beads of the former s-shell.
By iterating this, we get infinitely many aligned s-beads, all adjacent with the two 1-beads and two r-beads of the inital s-shell.
This is incompatible with $s>0$.
This last case is thus, as the three previous ones, excluded.

\paragraph{Dimension one: remaining cases.}
For the other $10$ cases, the combinatorics of the 1s- and rs-necklaces characterize the s-shell: 
\begin{itemize}
\item a tetrahedron with r- and 1-beads (not all identical):
$$
rrr/1rr,
\qquad
1rr/11r,
\qquad
11r/111.
$$
\item bipyramids with a base of 3 to 5 r-beads and two 1-beads on apexes:
$$
rrr/1r1r,
\qquad
rrrr/1r1r,
\qquad
rrrrr/1r1r.
$$
\item bipyramids with a base of 3 to 5 1-beads and two r-beads on apexes:
$$
1r1r/111,
\qquad
1r1r/1111,
\qquad
1r1r/11111.
$$
\item
an octahedron with three mutually adjacent r-beads and three 1-beads:
$$
11rr/11rr
$$
\end{itemize} 

For the tetrahedra of r- and 1-beads, the only way to tile the space with such tetrahedra is to form a close-packing of 1-beads with r-beads filling the octahedral holes.
Indeed, otherwise we would get another compact packing by two sizes of spheres - this is impossible \cite{Fer19}.
In this packing, the tetrahedra created by the r-beads which fill octahedral holes are 111r (case 11r/111).
Replacing $r$ by its value $r=\sqrt{2}-1$ in the equations associated with the necklaces show that $s$ is the smallest positive roots of $2s^4 + 4s^3 - 14s^2 + 8s - 1$:
$$
s=\tfrac{1}{2}\sqrt{9+6\sqrt{2}}-\sqrt{2}-\tfrac{1}{2}\simeq 0.18
$$
This yields the third case of Theorem~\ref{th:main}.\\

The bipyramid symmetry ensures that the s-bead is in the base plane.
Its size is thus determined by the base beads.
Namely, for a 1-bead base of size $k$:
$$
s=s_k:=\frac{1-\sin\tfrac{\pi}{k}}{\sin\tfrac{\pi}{k}}.
$$
For base of r-beads, this yields the same values for $\tfrac{s}{r}$.
One can then use the equation associated with 1s- or rs-necklace to compute $r$.
Only two cases yield values $r$ and $s$ such that $0<s<r<1$:
\begin{enumerate}
\item
The bipyramid with a base of three 1-beads yields
$$
r=\tfrac{1-\sqrt{3}}{2}\simeq 0.21
\quad\textrm{and}\quad
s=s_3=\tfrac{2}{3}\sqrt{3}-1\simeq 0.15
$$
\item
The bipyramid with a base of five r-beads yields for $r$ the smallest root of $r^4 - 20r^3 + 45r^2 - 30r + 5$ and $\tfrac{s}{r}=s_5$:
$$
r=5-2\sqrt{5} - \tfrac{1}{2}\sqrt{130-58\sqrt{5}}\simeq 0.25
\quad\textrm{and}\quad
s=\tfrac{1}{2}\sqrt{10-2\sqrt{5}}-1\simeq 0.17
$$
\end{enumerate}
Unfortunately, in both cases, an exhaustive search (using interval arithmetic and the above values of $r$ and $s$) shows that there are no r-shell, hence no compact packing with all the three sizes of spheres (maybe these s-shells can however be used to form packings with high density, though non-compact).\\

For the octahedron 11rr/11rr, one has $s\geq s_3\simeq 0.15$, otherwise the s-bead cannot be in contact with the three mutually adjacent 1-beads of the s-shell.
This lower bound yields an upper bound on the maximal number of beads in any necklace.
Namely there are at most $6$ beads.
An exhaustive search yields $162$ necklaces with $3$ to $6$ beads.
Since the octahedral s-shell contains two adjacent r-beads, there must be an rr-necklace.
An exhaustive search shows that the only possible rr-necklace is 1srrs. 
The associated equation, together with the ones associated with the 1s- and rs-necklaces, form a zero-dimensional system.
This system has a unique solution $(r,s)$ such that $0<s<r<1$: $r$ is the smallest root of $x^2-8x+2$ and $s$ is the smallest root of $x^4-24x^3+96x^2-96x+16$:
$$
r=4-\sqrt{14}\simeq 0.26
\quad\textrm{and}\quad
s=6-3\sqrt{2}-\sqrt{42-28\sqrt{2}}\simeq 0.21
$$
These values allow an exhaustive search for necklaces and shells.
For necklaces:
\begin{center}
\begin{tabular}{c|cccccc}
body/head & 1s & rs & ss & rr & 1r & 11\\
\hline
necklaces & 11rr & 11rr & $\emptyset$ & 1srrs & 11srs & 111r1srs
\end{tabular}
\end{center}
The above octahedral s-shell is the unique s-shell.
There is a unique r-shell: it has $12$ beads, three of each size, centered on the vertices of the unique polyhedron which has the following $14$ faces:
\begin{center}
\begin{tabular}{c|ccccc}
face & 111 & rrr & 11s & rrs & 1rs\\
\hline
number & 1 & 1 & 3 & 3 & 6
\end{tabular}
\end{center}
The unique r-shell has an rrr triangular face which, together with the center of the shell, form a tetrahedra.
The r-beads thus always form tetrahedra.
Consider such a tetrahedron.
The unique r-shell structure ensures that each face of this tetrahedron is a face of an s-shell.
The 1-beads of these four s-shells then form a regular octahedron.
Four (non-adjacent) faces of this octahedron are the 111 triangular faces of the four s-shells.
Each of the four other faces forms, together with an r-bead of the central tetrahedron (namely the r-bead shared by the three s-shells which contain the 1-beads of the considered face), a tetrahedron 111r.
Any compact packing by spheres of these sizes can thus be seen as a tiling by regular octahedra (the previous one, which contains a tetrahedron of r-beads and four s-beads) and, as usual, regular tetrahedra.
This is thus a close-packing of unit spheres.
This yields the last case of Theorem~\ref{th:main}.\\

The exhaustive search for 1-shells seems to be challenging (we stopped it after $48$ hours of computation).
Nevertheless, since we know the possible packings, we can now describe them.
There are two 1-shells.
Both contain $36$ beads, $12$ of each size, centered on the vertices of a polyhedron which has $68$ faces.
The two polyhedra are obtained by modifying in the same way the $6$ square faces of either a cuboactahedron or a triangular orthobicupola, both with 1-beads centered on their vertices.
Namely, two r-beads are inserted in the square of 1-beads to form two 11r triangle and two trapeze 11rr.
An s-bead is inserted in each trapeze to form two 1rs triangles, one 11s triangle and one rrs triangle.

\bibliographystyle{alpha}
\bibliography{3balls}

\end{document}